\documentclass[12pt]{article}
\usepackage[utf8]{inputenc}
\usepackage[english]{babel}
\usepackage{amsfonts}
\usepackage[titletoc]{appendix}
\usepackage{amssymb}

\usepackage{amsmath}
\usepackage{amsthm}
 \usepackage{CJK}
\numberwithin{equation}{section}

\usepackage{graphicx}
\newcommand{\simeqd}{\mathrel{\rotatebox[origin=c]{-90}{$\simeq$}}}

\begin{document}

 \title{Rational curves on complete intersection Calabi-Yau 3-folds }
\author{ B. Wang \\
\begin{CJK}{UTF8}{gbsn}
(汪      镔)
\end{CJK}}
\date {}

\newcommand{\Addresses}{{
  \bigskip
  \footnotesize

   \textsc{Department of Mathematics, Rhode Island college, Providence, 
   RI 02908}\par
  \text{E-mail address}:  \texttt{binwang64319@gmail.com}

}}

\renewcommand{\thefootnote}{\arabic{footnote}}

\newtheorem{thm}{Theorem}[section]

\newtheorem{ass}[thm]{\bf {Assumption} }
\newtheorem{proposition}[thm]{\bf {Proposition} }
\newtheorem{theorem}[thm]{\bf {Theorem} }
\newtheorem{fo}[thm]{\bf Focus}
\newtheorem{corollary}[thm]{\bf Corollary}
\newtheorem{definition}[thm]{\bf Definition}
\newtheorem{lemma}[thm]{\bf Lemma}

\newcommand{\xdownarrow}[1]{%
  {\left\downarrow\vbox to #1{}\right.\kern-\nulldelimiterspace}
}

\newcommand{\xuparrow}[1]{%
  {\left\uparrow\vbox to #1{}\right.\kern-\nulldelimiterspace}
}

\bigskip

\maketitle

\begin{abstract}   
We prove the following results. 
If $X_3$ is a generic complete intersection  Calabi-Yau 3-fold, 
\par

(1) then for each natural number $d$
 there exists a rational map \par\hspace{1 cc} $c\in Hom_{bir}(\mathbf P^1, X_3)$ of $deg(c(\mathbf P^1))=d$, 
 \par
(2)  further more all such $c$ are immersions satisfying
\begin{equation}
N_{c(\mathbf P^1)/ X_3}\simeq \mathcal O_{\mathbf P^1}(-1)\oplus  \mathcal O_{\mathbf P^1}(-1).
\end{equation}

 \bigskip

 \par

\end{abstract}

\tableofcontents

\bigskip

\section{Statement}

\bigskip

 In mirror symmetry there is a general consensus  that a ``generic'' Calabi-Yau 3-fold over $\mathbb C$
 should contain and only contain finitely many  irreducible rational curves  of
each degree with respect to the polarization.  This is the basis for the calculation of instantons in Mirror symmetry.
   While we are still in the process of finalizing the  precise statement, let us prove it for the complete intersection Calabi-Yau 3-folds in 
a single projective space.
\bigskip

\begin{theorem}\quad\par 
Let $ X_3$ be a generic, complete intersection Calabi-Yau 3-folds over $\mathbb C$. \par
Then\par

(1) $X_3$ admits an irreducible rational curve $C$ of each degree,\par
(2) all such $C\subset  X_3$ are immersions and the normal bundle $N_{C/X_3}$ is  \par\hspace{1 cc}  isomorphic to 

\begin{equation} \mathcal O_{\mathbf P^1}(-1)\oplus \mathcal O_{\mathbf P^1}(-1).\end{equation}
\end{theorem}
\bigskip

\section{Sketch of the proof}

\subsection{Setting}
Throughout the paper rational curves are curves rationally parametrized by $\mathbf P^1$.  
The image is called an irreducible rational curve. \medskip

Rational curves on projective varieties has been a topic  for many decades. The general theory which
has its own technique  is not the focus of this paper.   Instead we are interested in a specific type of problems with a quite different technique.
\par
\begin{fo} Which  generic complete intersection $X$ admits an irreducible rational rational curve $C$ of each degree?
\end{fo}\bigskip

\begin{fo}   Once the first question is affirmative, what is the normal sheaf
\begin{equation}
N_{C/X}?
\end{equation}
\end{fo}

Theorem 1.1  tries to answer these two questions  in the case of Calabi-Yau 3-folds.

\bigskip

The idea of the work starts from and stays in  a down-to-earth setting, which employees  linear algebra only.   The method 
  first converts the invariant expression of Theorem 1.1 to
a variant expression as the content of Theorem 1.1 stays the same.  Then it explores the unique linear algebra in the variant setting
 to reach an algebraic result. At last it converts
the algebraic result back to the invariants.   In technique the first conversion 
\begin{equation}
{\text  Invariant \Rightarrow Variant }
\end{equation}
uses classical geometry. The second conversion
 \begin{equation}
{\text  Variant \Rightarrow invariant }
\end{equation}
uses Clemens' deformation idea [1]. 
\bigskip

Let's start with this alternative setting.  
 Let
\begin{equation}M_d=( H^0(\mathcal O_{\mathbf P^1}(d))^{\oplus (n+1)}\simeq \mathbb C^{(n+1)(d+1)}.\end{equation}
The open set $M_{bir, d}$ of $M_d$ represents (but is not equal to)  $$
\{  c\in Hom_{bir}(\mathbf P^1, \mathbf P^n): deg(c(\mathbf P^1))=d \}.$$
Let $$g_i, i=1, \cdots, r=n-3$$ be  sections in $H^0(\mathcal O_{\mathbf P^n}(h_i))$ and 
\begin{equation}  
\sum_{i=1}^r h_i =n+1. 
\end{equation}
Let
 \begin{equation}
X_i=\cap_{k=1}^{k=n-i} div(g_k).
\end{equation}
In this paper, the Cartesian product 
$$(g_1, \cdots, g_{n-i})$$ is also called 
a complete intersection of type $(h_1, \cdots, h_{n-i})$. 
So
\begin{equation}
X_3=\cap_{i=1}^r div(g_i)
\end{equation}
is a complete intersection Calabi-Yau 3-fold  in the usual sense for generic 
$$(g_1, \cdots, g_r)\in \prod_{i+1}^r H^0(\mathcal O_{\mathbf P^n}(h_i)). $$ 

\bigskip

  Choose  distinct  $h_id+1$ points  $t^i_j\in \mathbf P^1$.  Let
\begin{equation} \mathbf t^i=(t^i_1, \cdots, t^i_{h_id+1})\in Sym^{h_id+1}(\mathbf P^1).\end{equation}
In the rest of the paper,  we'll use following conventions in affine coordinates.
 \par 
(a) $t^i_j$ or $t$  denotes a  complex number  which is a point in an affine open \par\hspace{1cc}  set  $\mathbb C\subset \mathbf P^1$, \par (b)
 $c(t)$ denotes the image 
$$\begin{array}{ccc} \mathbb C &\stackrel{c}\rightarrow & \mathbb C^{(n+1)(d+1)}\\
t &\rightarrow & c(t),
\end{array}$$\par 
(c)  $g_i\in H^0(\mathcal O_{\mathbf P^n}(h_i))$ is a homogeneous polynomial of degree $h_i$ in \par\hspace{1cc} $n+1$ variables. 
\par

We should note that in these affine coordinates, the incidence relation $c^\ast(g_i)=0$ can be expressed
as the composition $g_i(c(t))=0$ for all $$t\in \mathbb C.$$
Let $C_M$ be a system of affine coordinates for $M_d$, which determines an isomorphism 
$$M_d\simeq \mathbb C^{(n+1)(d+1)}. $$

We define a system of polynomials in the variable $c\in M_d$
\begin{equation}
g_i(c(t^i_j)), i=1, \cdots, r, j=1, \cdots, h_id+1.
\end{equation}
Then the subsets of polynomials $$
g_1, \cdots, g_l, l\leq r$$
give a rise to a holomorphic map $\mu_l$
\begin{equation}\begin{array}{ccc}
\mu_l: M_d\simeq \mathbb C^{(n+1)(d+1)}\to \mathbb C^{m_ld+l},\end{array}\end{equation}
where $m_l=\sum_{i=1}^l h_i$. 
We''ll denote 
\begin{equation}
I_l=\mu_l^{-1}(0)
\end{equation}
which will be called the incidence scheme of rational curves 
on the complete intersection of
$g_1, \cdots, g_l$. Then the differential map $(\mu_l)_\ast$ is 
represented by the  Jacobian matrix of size 
$$\biggl(m_ld+l\biggr)\times \biggl(n+1\biggr)\biggl(d+1\biggr)$$ denoted by
$J_l$, which depends on $g_i, \mathbf t^j$.  However once the points $\mathbf t^j$ are fixed, 
the matrix $J_l$ is well-defined and varied algebraically   on the entire affine
space $M_d\times \mathbb A$, where $\mathbb A\subset  \prod_i \mathbf P(H^0(\mathcal O_{\mathbf P^n}(h_i)))$ is affine. 
We call it the 
Jacobian matrix of the incidence scheme $ I_l$.   
\par

The scheme $I_l$, which could be reducible with multiple dimensions,\footnote {For instance, it contains the components of 
multiple covering maps of $\mathbf P^1$.}   is the alternative to various moduli spaces 
of rational curves and maps.  
In this paper we show a methodology in the calculation of the Jacobian matrix $J_l$ at  a point of $I_l$ corresponding to
an irreducible rational curve on  generic complete intersections.  The components containing such points always have 
the smallest dimension as expected. They correspond to
the components that  have actual  fundamental classes  instead of  virtual fundamental classes  in Mirror symmetry.  
The methodology is rooted in a specific pattern of the Jacobian matrix $J_l$ modeled on the Vandermonde matrices. 
\par

In general we are interested in the following questions,\par
(1) what is the dimension of $I_l$?, \par
(2) is it reduced? \par
(3) is it irreducible, and if not, what is the structure of each component? \par

Answers to these questions would solve some technical problems in Mirror symmetry,  that have been left out of
 the physics'  formulation.  
Thus it is an alternative to patch some holes in the general mathematical  theory.   
The complete answers to these questions are out of scope of this paper.  
Here we are going to use linear algebra to explore the dimension of  $I_l$ for some complete intersections.
\bigskip

\subsection {Existence }
In the first part, we prove the existence of irreducible rational curves $C$ 
on a generic complete intersection Calabi-Yau 3-fold $X_3$. We avoid a direct construction.\bigskip

(I) First we'll use the proven existence of irreducible rational curves 
of arbitrary degrees on a single  generic hypersurface of lower degrees (which are Fano).
 After its extension to a special complete intersection 
in  $\mathbf P^n$  by joining more Fano hypersurfaces, we use linear algebra to glue
the block matrix for each Fano hypersurface to obtain the non-degeneracy of Jacobian matrix $J_r$. 
This will show the existence of an irreducible rational curve $C'$ of each degree 
on a special complete intersection Calabi-Yau 3-fold $X_3'$ 
with the non-degenerate  Jacobian $J_r'$ at the point $(C', X_3')$.
\footnote{ For curves of  a fixed and small degree, 
computer software  has been used to find its Jacobian. See p. 295, [2].} 
\par

(II) Applying Clemens' deformation idea, this smooth $C'$ is deformed to generic 
$X_3$ as a different irreducible rational curve $C\subset X_3$. The requirement  for such a deformation is the non-degeneracy of 
Jacobian matrix $J_l$. \bigskip

\subsection{Rigidity}
The second step is to show the normal bundle of all such irreducible curves $C\subset X_3$ are split  as in (0.1).
This rigidity is determined by the dimension of the incidence scheme $I_r$. 
First we have the general study of a uniruled projective variety   to deduce  that 
a free rational curve on it has the unobstructed deformation on its generic hypersurfaces. This no-trivial result follows 
from [7] or [8].   
Then we notice that the filtration of complete intersections
\begin{equation}
X_3\subset \cdots \subset X_0\simeq \mathbf P^n
\end{equation}
corresponds to another filtration of subvarieties
\begin{equation}
I_r|_{c}\subset I_{r-1}|_{c}\subset \cdots\subset I_1|_{c}\subset I_0|_{c},
\end{equation}where the subscript $|_{c}$ means
the analytic neighborhood around the point $c$ with $c(\mathbf P^1)=C$.   
Then we notice all $X_i$ for $i\neq 3$ are Fano, therefore uniruled. 
Applying above deformation result for rational curve on uniruled varieties, we obtain a recursive formula
\begin{equation}
dim( I_{l}|_{c})= dim(I_{l+1}|_{c})-(h_ld+1).
\end{equation}
( or equivalently $dim(H^0(c^\ast (T_{X_l}))=dim(H^0(c^\ast (T_{X_{l+1}}))-(h_ld+1).$). 
Applying the Calabi-Yau condition (2.5)  we obtain that
\begin{equation}
dim(I_r|_{c})=4.
\end{equation}
(or $dim( H^0(c^\ast (T_{X_3}))=3$). 
\bigskip

Then the dimension $dim( H^0(c^\ast (T_{X_3}))=3$ forces $c$ to be an immersion and furthermore 
to be rigid. 

\bigskip

We organize the rest of the paper as follows.  In section 3, we prove the existence, and in section 4 we prove the rigidity. 
Appendix covers  a particular technique 
in linear algebra for the existence. 

\section{Existence of rational curves}

\bigskip

Starting from this section we give technical proofs. This section focuses on the existence. \bigskip

We resume all notations introduced in section 1. Recall
 $g_i, i=1, \cdots, r$ are  sections in $H^0(\mathcal O_{\mathbf P^n}(h_i))$ such that
$n-3=r$ and $n+1=\sum_i h_i$.  Let
\begin{equation}
X_l=\cap_{i=1}^{n-l} div(g_i), l\leq r
\end{equation}
 In particular, if 
$$(g_1, \cdots, g_r)\in \prod_{i=1}^r H^0(\mathcal O_{\mathbf P^n}(h_i))$$ is generic, $X_3$
is a smooth complete intersection Calabi-Yau 3-fold.

In this subsection, we prove part (1) of Theorem 1.1. It asserts the existence of 
an irreducible rational curve of each degree $d$ on the generic $X_3$. 

\bigskip

\begin{theorem}
There exist irreducible rational curves of each degree $d$ on the following generic 
complete intersection Calabi-Yau 3-folds.
\begin{equation}
(2, 2, 2,2) \ type \ in \ \mathbf P^7
\end{equation}
\begin{equation}
(3, 2,2)  \ type \ in \ \mathbf P^6, 
\end{equation}
\begin{equation}
(3, 3) \ and  \ (4, 2)   \ types \ in \ \mathbf P^5 \end{equation}
\begin{equation}
(5) \ type, i.e.  quintic\ 3\text {-}fold \ in \ \mathbf P^4.\end{equation}

\end{theorem}

\bigskip

The case (3.5) has been proved by H. Clemens [1] and S. Katz ([3]). Let's consider other 3 cases. 
Suppose that
in all 3 cases, there exist  a special complete intersection, denoted by sections 
$g_1', \cdots, g_r'$  and 
a smooth rational curve $c_g$ of degree $d$ with $c_g^\ast(g_i')=0$ for all $i$ such that the corresponding Jacobian matrix 
$J_r$ at $c_g$ has full rank.  Then we divide the coordinates of $M_d$ to two parts, $c_{ind}$ and $c_{free}$.
The partial derivatives with respect to $c_{ind}$ form the maximal minor block $J_{ind}$ in the Jacobian matrix $J_r$.
Thus  $J_{ind}$ is a maximal non-degenerate block.   The size of $J_{ind}$ is therefore $$ \biggl ((n+1)d+r\biggr) \times  \biggl ((n+1)d+r\biggr)$$
Let the remaining coordinates of $M_d$ be $c_{free}$, and corresponding block matrix in $J_r$ is $J_{free}$.
 Hence $(c_{ind}c_{free})$ are the affine coordinates of $M_d$ and $(J_{ind}J_{free})=J_r$. 
Since  $J_{ind}$ evaluated at the special complete intersection $(g_i')$ and smooth $c_g$ 
 is non-degenerate,   by the implicit function theorem in complex analysis, there exists
analytic functions near $c_g$, $c_{ind}=\alpha_{ind}( c_{free}, g_1, \cdots, g_r)$ ($g_i$ are locally free sections) such that
\begin{equation}
 g_i(c(t^i_j))=0, \ for\ all\ i, j
\end{equation}
where $c=\alpha_{ind}( c_{free}, g_1, \cdots, g_r) c_{free}$.
In geometric term, this means that 
there exists a smooth rational curve $c$ of degree $d$ on each same type of complete intersection (free choice of $g_i$)
 near the special one $(g_i')$. Since these deformed complete intersections $(g_i)$ are all generic in Zariski-topology, we complete the proof of the existence in these cases. So in the following subsections  we are going to find such a special complete intersection in each case.

 \bigskip

 \bigskip

\subsection{Complete intersection of $(2, 2, 2, 2)$ type in $ \mathbf P^7$}

\begin{proposition}

There is a smooth rational curve $C_g$ of each degree $d$ on a special  complete
intersection $X_3$ of type $(2, 2,2,2)$ in $\mathbf P^7$ such that $J_4|_{c_g}$ at $(c_g, X_3)$ is non-degenerate, 
where $c_g$ is the normalization of $C_g$.    

\end{proposition}
\bigskip

\begin{proof} Let $[z_0, \cdots, z_7]$ be homogeneous coordinates of $\mathbf P^7$. Let 
$\mathbf P^3$ the subspace defined by $z_4=z_5=z_6=z_7=0$. 
First we consider a generic quadric $g\in H^0(\mathcal O_{\mathbf P^3}(2))$ in $\mathbf P^3$.
By [4], Hilbert scheme $Hilb_{c(\mathbf P^1)}(div(g))$ at a smooth rational curve 
$$c\in Hom_{bir}(\mathbf P^1, \mathbf P^n)$$ of degree 
$d$ is non-empty and smooth with the expected dimension. 
 Then $I_1$ is smooth at $c$ with expected dimension. Let $J_1'$ be the corresponding
Jacobian matrix of $I_{1}$ in $\mathbf P^3$ (as defined in the introduction). By Lemma A.2, it has full rank. 
In the following we use the ``gluing" technique for block matrices to extend the Jacobian matrix $J_1'$ to  $\mathbf P^7$.
 Assume $c=[c_0, c_1, c_2, c_3]$. It is extended to a smooth rational curve $c_g$ in $\mathbf P^7$ as follows
\begin{equation}
c_g=[c_0, c_1, c_2, c_3, c_0, c_1, c_2, c_3],
\end{equation}
which is isomorphic to $c$.  We define the special complete intersection of $(2, 2, 2, 2)$ type as follows.
Let \begin{equation}\begin{array}{c}
g_1 =g (by\ the \ extension),\\ 
g_2=g(z_4, z_5, z_2, z_3)\\
 g_3=g(z_0, z_1, z_6, z_7), \\
g_4=z_6^2+z_7^2-z_2^2-z_3^2.\end{array}\end{equation}
be four quadrics in $\mathbf P^7$. Then we have  $c_g^\ast (g_i)=0$ for all $i$.

Let $\theta_1, \theta_3, \theta_5, \theta_7$ be the affine coordinates  for $2nd, 4th, 6th$ and $8th$ 
copies $H^0(\mathcal O_{\mathbf P^1}(d))$ in $M_d$. 
Let $V_0$ be the hyperplane of $H^0(\mathcal O_{\mathbf P^1}(d))$, whose elements vanish at $0$. Let
 $\theta_0, \theta_2, \theta_4, \theta_6$ be the affine coordinates for $V_0$ in the $1st, 3rd, 5th, 7th$ components of $M_d$.
 The we'll show these variables,  that are not all variables of $M_d$, are $c_{ind}$. 
This amounts to show that the Jacobian matrix has full rank. Thus we first write down the blocks of Jacobian matrix. 
We let \begin{equation}\begin{array}{c}
B_{11}={\partial (g_1(c_g(t^1_1)), \cdots, g_1(c_g(t^1_{2d+1})))\over \partial (\theta_0, \theta_1)}\\
B_{12}={\partial (g_1(c_g(t^1_1)), \cdots, g_1(c_g(t^1_{2d+1})))\over \partial (\theta_2, \theta_3)}\\
B_{22}={\partial (g_2(c_g(t^2_1)), \cdots, g_1(c_g(t^2_{2d+1})))\over \partial (\theta_2, \theta_3)}\\
B_{23}={\partial (g_2(c_g(t^2_1)), \cdots, g_2(c_g(t^2_{2d+1})))\over \partial (\theta_4, \theta_5)}\\
B_{33}={\partial (g_3(c_g(t^3_1)), \cdots, g_3(c_g(t^3_{2d+1})))\over \partial (\theta_4, \theta_6)}\\
B_{34}={\partial (g_3(c_g(t^3_1)), \cdots, g_3(c_g(t^3_{2d+1})))\over \partial (\theta_6, \theta_7)}\\
B_{42}={\partial (g_4(c_g(t^4_1)), \cdots, g_g(c_g(t^4_{2d+1})))\over \partial (\theta_2, \theta_3)}\\
B_{44}={\partial (g_4(c_g(t^4_1)), \cdots, g_4(c_g(t^4_{2d+1})))\over \partial (\theta_6, \theta_7)}

\end{array}\end{equation} 
Then one of  maximal minor blocks  $J_{ind}$ of the Jacobian matrix $J_4$ ( in $\mathbf P^7$)  at $c_g$ is formed by the block matrices
\begin{equation}J_{ind}=\left (\begin{array}{cccc} 
B_{11} &B_{12} &0 & 0 \\
0&B_{22}& B_{23} &0  \\
0 &0& B_{33}  & B_{34}\\
0 &B_{42}&0  & B_{44}\\
\end{array}\right)\end{equation} 
We can verify  that $B_{ij} $ in (3.10) satisfy all conditions in Lemma A.2.  Furthermore  all conditions in  
Lemma A.6 for $J|_{c_g}$ are also satisfied.  By the Lemmas $J|_{c_g}$  is non-degenerate.  
 We complete the proof. 
\end{proof}

\bigskip

\subsection{Complete intersection of $(3, 2, 2)$ type in  $\mathbf P^6$}

\begin{proposition}

There is a smooth rational curve $C_g$ of each degree $d$ on a special  complete
intersection $X_3$ of type $(3, 2,2)$ in $\mathbf P^6$ such that $J_3|_{c_g}$ at $(c_g, X_3)$ is non-degenerate, 
where $c_g$ is the normalization of $C_g$.    

\end{proposition}
\bigskip

\begin{proof}

Let $z_0, \cdots, z_6$ be homogeneous coordinates of $\mathbf P^6$. 
Let $\mathbf P^3\subset \mathbf P^6$ be the subspace defined by $z_4=z_5=z_6=0$. 
  By [4], there exists a generic $$g\in H^0(\mathcal O_{\mathbf P^3}(2))$$ such that it
contains a smooth rational curve 
\begin{equation}
c=[c_0, c_1, c_2, c_3]. 
\end{equation}

We define a new smooth rational curve of degree $d$ to 
be 
\begin{equation}
c_g=[c_0, c_1, c_2, c_3, c_0, c_1, c_2]
.\end{equation}

We define hypersurfaces  as
\begin{equation}\begin{array}{c}
g_1=z_0g(z_0, z_1, z_2, z_3)+(z_0-z_4)\alpha(z_0, z_1, z_2, z_3)\\
g_2=g(z_4, z_5, z_2, z_3)\\
g_2=g(z_4, z_5, z_6, z_3)
\end{array}\end{equation}
where $\alpha$ is a generic quadric with respect to $g$. 
Let $\theta_1,\cdots, \theta_7$ be the affine coordinates for $7$ 
copies $H^0(\mathcal O_{\mathbf P^1}(d))$ in $M_d$. 
Let \begin{equation}\begin{array}{c}
\mathbf t^1=(t_1^1, \cdots, t^1_{3d+1})\in sym^{3d+1}(\mathbf P^1)\\
\mathbf t^2=(t_1^2, \cdots, t^2_{2d+1})\in sym^{2d+1}(\mathbf P^1)\\
\mathbf t^3=(t_1^3, \cdots, t^3_{2d+1})\in sym^{2d+1}(\mathbf P^1).
\end{array}\end{equation}

We write down the block matrices for the Jacobian matrix. 
Let
\begin{equation}\begin{array}{c}
B_{1j}={\partial (g_1(c_g(t^1_1)), \cdots, g_1(c_g(t^1_{3d+1})))\over \partial \theta_j}\\
B_{2j}={\partial (g_2(c_g(t^2_1)), \cdots, g_2(c_g(t^2_{2d+1})))\over \partial \theta_j}\\
B_{3j}={\partial (g_3(c_g(t^3_1)), \cdots, g_3(c_g(t^3_{2d+1})))\over \partial \theta_j}.\end{array}\end{equation}

Next we write down the Jacobian matrix directly as 
 
\begin{equation}J_3=\left (\begin{array}{ccccccc} 
B_{11} &  B_{12} &  B_{13} & B_{14} & B_{15} & 0 & 0\\
0 &  0 & B_{23} & B_{24} & B_{25} & B_{26} & 0\\
0 &  0 & 0 & B_{34} & B_{35} & B_{36} & B_{37}
\end{array}\right)\end{equation}
Then we verify $B_{ij}$ satisfy all conditions in Lemmas A.2, A.3. Then we can apply  
 Lemma A.5.  
We obtain that  Jacobian matrix $J_3$ of the incidence scheme has full rank. We complete the proof of this case.
\end{proof}

\bigskip

\subsection{Complete intersections of $(3,3)$ and $(4, 2)$  types in  $\mathbf P^5$}
\bigskip

\begin{proposition}

There is a smooth rational curve $C_g$ of each degree $d$ on a special  complete
intersection $X_3$ of type $(3, 3)$ in $\mathbf P^5$ such that $J_2|_{c_g}$ at $(c_g, X_3)$ is non-degenerate, 
where $c_g$ is the normalization of $C_g$.    

\end{proposition}
\bigskip

\begin{proof} Let $z_0, \cdots, z_5$ be homogeneous coordinates of $\mathbf P^5$. 
Let $\mathbf P^3\subset \mathbf P^5$ be the subspace defined by $z_4=z_5=0$. 
  By [4], there exists a generic $g\in H^0(\mathcal O_{\mathbf P^3}(2))$ such that it
contains an irreducible rational curve 
\begin{equation}
c=[c_0, c_1, c_2, c_3]
\end{equation}
of degree $d$.
Next we define a smooth rational curve of degree $d$ in $\mathbf P^5$ as
\begin{equation}
c_g=[  c_0, c_1, c_2, c_3, c_1, c_0].
\end{equation}
We also define two cubic hypersurfaces
\begin{equation}\begin{array}{c}
g_1=z_0 g(z_0, z_1, z_2, z_3)+(z_0-z_5)\alpha(z_0, z_1, z_2, z_3) \\
g_2=z_0 g(z_5, z_4, z_2, z_3)+(z_1-z_4)\alpha(z_0, z_1, z_2, z_3),
\end{array}\end{equation}
where $\alpha$ is a generic quadric with respect to $g$. 
 Let $\theta_1,\cdots, \theta_6$ be the affine coordinates  for $6$ 
copies $H^0(\mathcal O_{\mathbf P^1}(d))$ in $M_d$ (all variables of $M_d$). 
Let \begin{equation}\begin{array}{c}
\mathbf t^1=(t_1^1, \cdots, t^1_{3d+1})\in sym^{3d+1}(\mathbf P^1)\\
\mathbf t^2=(t_1^2, \cdots, t^2_{3d+1})\in sym^{3d+1}(\mathbf P^1)
.\end{array}\end{equation}
Let 
\begin{equation}
B_{ij}={\partial (g_i(c_g(t^i_1)), \cdots, g_i(c_g(t^i_{3d+1})))\over \partial \theta_j}
\end{equation}
for $i=1, 2$, $j=1, \cdots, 6$.  
Then we directly write down the Jacobian $J_2$ of the incidence scheme
at $c_g$. \begin{equation}J_2=\left (\begin{array}{cccccc} 
B_{11} &  B_{12} &  B_{13} & B_{14} &0 & B_{16}\\
0 & B_{22} & B_{23} & B_{24} & B_{25} & B_{26}
\end{array}\right)\end{equation}
Similarly  conditions of Lemma A.3 are met. 
By Lemma A.4, it has full rank.  We complete the proof in this case.
\end{proof}

\bigskip

\begin{proposition}

There is a smooth rational curve $C_g$ of each degree $d$ on a special  complete
intersection $X_3$ of type $(4, 2)$ in $\mathbf P^5$ such that $J_2|_{c_g}$ at $(c_g, X_3)$ is non-degenerate, 
where $c_g$ is the normalization of $C_g$.    

\end{proposition}
\bigskip

  \begin{proof}  As in the case 1, 
$z_0, \cdots, z_5$ are homogeneous coordinates of $\mathbf P^5$. Let 
$\mathbf P^4$ be the subspace covered by the coordinates
$$[z_0, z_1, z_2, z_3, z_4, 0].$$
Let $\mathbf P^3$ be the subspace covered by coordinates
$$[z_0, z_1, z_2, z_3, 0, 0].$$
By the Mori's result [6]. there is a smooth rational curve $c_s$ of degree $d$ on a generic quartic 
$g\in H^0(\mathcal O_{\mathbf P^3}(4))$.  Let \begin{equation}
\mathbf t^1=(t_1^1, \cdots,  t_{4d}^1)\in sym^{4d}(\mathbf P^1)\end{equation}
be generic. 

Then the Jacobian matrix $J'$
\begin{equation}
J'={\partial \bigl( g(c_s(t_1^1)), \cdots, g(c_s(t_{4d}^1))\bigr)
\over \partial (\theta_0, \cdots, \theta_3)}
\end{equation}
has full rank. \footnote{ The number of points $t^1_j$ is unusual. 
This exceptional case holds only for quartic surface in $\mathbf P^3$.} Let's extend it to $\mathbf P^4$. 
Define a new quartic  by setting $g_1=g+z_0^3 z_4$ in $\mathbf P^4$ and new rational curve by setting  $$c=[c_0, c_1, c_2, c_3, 0]$$ 
where $[c_0, c_1, c_2, c_3]$ is the Mori's curve.
Next we add one generic point $$t^1_{4d+1}\in \mathbf P^1,$$ and
extend the coordinates of $$H^0(\mathcal O_{\mathbf P^1}(d))^{\oplus 4}$$ to $$H^0(\mathcal O_{\mathbf P^1}(d))^{\oplus 5}.$$ 
These add the a new row --the differential of $\mathbf d g_1(c(t_{4d+1}^1))$ to the Jacobian $J'$,  
Hence we  obtain a new Jacobian matrix $J_1$ in $\mathbf P^4$ at $c$, 
\begin{equation}J_1|_{c_g}=\left (\begin{array}{cc} 
J' &B_{12}\\
 B_{21}& B_{22} 
\end{array}\right)\end{equation} 
where
\begin{equation}\left (\begin{array}{cc} 
B_{21}& B_{22} 
\end{array}\right)\end{equation} 
 is the differential 1-form $\mathbf d g_1(c(t_{4d+1}^1))$. 
Hence $J_1|_{c}$ in $\mathbf P^4$ has full rank. 
 To summarize it, we found a special quartic $g_1$ in $\mathbf P^4$ containing a smooth 
rational curve $c$ of the given degree and its Jacobian matrix is of full rank. \footnote{ We added
one more feature to Mori's existence. That is the non-degeneracy of the Jacobian matrix. }
\par

Next  we extend it one more time to $\mathbf P^5$. 
 Let $g_1$ be the extension of original $g_1$ to $\mathbf P^5$. 
Let $g_2=z_0z_4+z_1z_5$ and the $$c_g=[c_0, c_1, c_2, c_3, 0, 0]$$ be the smooth curve in $\mathbf P^5$.
Thus $c_g$ lies on the complete intersection of $g_1, g_2$. 
We add new set of generic $2d+1$ points in $\mathbf P^1$, 
\begin{equation}
\mathbf t^2=(t_1^2, \cdots, t_1^{2d+1})\in sym^{2d+1}(\mathbf P^1).
\end{equation}

Then we can write down the Jacobian matrix $J_2$ in this case. It is equal to
\begin{equation} J_2= \left (\begin{array}{cc} 
J_1 & A_{12}  \\
A_{21}& A_{22} 
\end{array}\right)\end{equation} 
where  \begin{equation}  \left (\begin{array}{cc} 
A_{21}& A_{22} 
\end{array}\right)\end{equation}  is the Jacobian matrix at $c_g$ of 
$2d+1$ many functions (in $c$), 
\begin{equation}
g_2(c(t_1^2)), \cdots, g_2(c(t_{2d+1}^2))
\end{equation}
and  $A_{22}$ is the block with respect to $2d+2$ coordinates in 
last two components of 
\begin{equation}
H^0(\mathcal O_{\mathbf P^1}(d))^{\oplus 6}.\end{equation}
By Lemma A.2, $ A_{22}$ has full rank. Hence the Jacobian matrix 
$J_2$ also has full rank. 
We complete the proof.\end{proof}
\bigskip

\subsection{Quintic in  $\mathbf P^4$}

\bigskip

This has been proved in [3]. We'll repeat the same construction in [1], but continue with our method.
 By the Mori's result [6]. there is a smooth rational curve $c_s$ of degree $d$ on a generic quartic $div(g)$
in $\mathbf P^3$.  Now we extend it to $\mathbf P^4$ by setting
\begin{equation}
c_g=[c_0, c_1, c_2, c_3, 0], for \ c_s=[c_0, c_1, c_2, c_3],
\end{equation} and the new 
quintic \begin{equation}
g_1=l g+q z_4
,\end{equation}
where $l$ is a generic linear polynomial and $q$ is a generic quartic. By [9], 
the corresponding Jacobian matrix $J_1$ has full rank.
\bigskip

\begin{corollary}
Let $X_3$ be a generic complete intersection Calabi-Yau 3-fold. Then it contains 
an irreducible rational curve of each degree $d$.

\end{corollary}

\bigskip

\begin{proof}
$X_3$ is projectively embedded as one of the complete intersections in Theorem 3.1. 
Therefore the corollary follows. 

\end{proof}

\bigskip

{\bf Remark} \quad
In above arguments for the existence, all rational curves are smooth. But non-smooth and irreducible rational curves
on a  generic complete intersection Calabi-Yau 3-fold do exist.

\bigskip

\section{ Rigidity of rational curves}

\subsection{Rational curves with unobstructed deformation}

In this subsection we let $Y$ be an arbitrary smooth 
projective variety over $\mathbb C$.  We'll prove a general assertion in Theorem 4.3 about uniruled variety.
For a parametrized rational curves on $Y$, 
the notion of having unobstructed deformation is weaker than being a 
 free morphism. Precisely let $c\in Hom_{bir}(\mathbf P^1, Y)$. \bigskip
 
 \begin{definition}

 If the normal 
sheaf $c^\ast (N_{c(\mathbf P^1)/Y})$ denoted simply by $N_{c/Y}$ has the vanishing first cohomology, i.e.  \begin{equation}
H^1(N_{c/Y})=0, \end{equation}
then we say $c$ has unobstructed deformation on $Y$. 
\end{definition}

\bigskip

Since $H^1(N_{c/Y})=0$
 is equivalent to $H^1(c^\ast(T_{Y}))=0$, then Definition 4.1  is  equivalent 
 to the following splitting
\begin{equation}
c^\ast(T_{Y})\simeq \oplus_j \mathcal O_{\mathbf P^1}(a_j), a_j\geq -1.
\end{equation}
On the other hand, we define

\bigskip

\begin{definition}
If $c^\ast(T_{Y})$ is generated by global sections, we say $c$ is a free morphism, and
$Y$ is uniruled. 
\end{definition}

\bigskip

This is a special case of the more general definition in [5].\bigskip

It is clear that $c$ being free is equivalent to the splitting 
\begin{equation}
c^\ast(T_{Y})\simeq \oplus_j \mathcal O_{\mathbf P^1}(a_j), a_j\geq 0.
\end{equation}
 So  being free is stronger than having unobstructed deformation.

\bigskip

\bigskip

\bigskip

\bigskip

\begin{theorem}  Let $\mathcal L\simeq \mathcal O_{\mathbf P^n}(1)|_Y$ be a very ample line bundle on $Y$, and $dim(Y)\geq 4$. 
Let $$X=div(f)\subset Y$$
where $f\in H^0(\mathcal L^h)$ is generic with an $h$.  
  Let 
$$ c: \mathbf P^1 \to C\subset  X\subset Y$$ be a birational morphism onto
an irreducible rational curve $C$.  If $c$ is free on $Y$, then $c$ has unobstructed deformation on $X$, i.e. 
\begin{equation}
  H^1(N_{c/X})=0.
\end{equation}
\end{theorem}

\bigskip

\begin{proof} of Theorem 4.3:
By the assumption of the theorem, we have a polarization of $Y$ such that 
\begin{equation}
Y\subset \mathbf P^n
\end{equation} is a smooth subvariety of dimension $\geq 4$, and
$\mathcal L=\mathcal O_{\mathbf P^n}(1)|_{Y}$. Let 
\begin{equation}  s\in H^0(\mathcal O_{\mathbf P^n}(h))
\end{equation}
be generic. 
Let $ f=s|_Y\in H^0(\mathcal L^h)$.

We denote \begin{equation}
div(f)=X, div(s)=Z.
\end{equation}
By the genericity of $Z$,  scheme-theoretically
\begin{equation}
X=Y\cap Z.
\end{equation}
Let $c:\mathbf P^1\to X$ be  generic  in $Hom_{bir}(\mathbf P^1, X)$.

We have a non-commutative diagram of exact sequences
\begin{equation}\begin{array}{ccccccccccc}
&&&&0 && 0&&&&\\
&&&&\downarrow &&\downarrow &&&&\\
0 &\rightarrow & H^0(c^\ast T_{X})&\rightarrow & H^0(c^\ast T_{Y})&\stackrel{\mu_1}\rightarrow &
H^0(c^\ast N_{X/Y})&\stackrel{\mu_2}\rightarrow &H^1 (c^\ast T_{X})
&\rightarrow &
H^1(c^\ast T_{Y})\\
&&&&\downarrow &&\simeqd &&&&\\
 0&\rightarrow & H^0(c^\ast T_{Z}) &\rightarrow & H^0(c^\ast  (T_{\mathbf P^n})) &\stackrel{\mu_3}\rightarrow & 
H^0(c^\ast (N_{Z/\mathbf P^n}))&\rightarrow &  H^1(c^\ast(T_{Z}) )&&\\
&&&& \downarrow   & & \|&  &&&\\
&&&& H^0(c^\ast (N_{Y/\mathbf P^n}))  &\stackrel{\mu_4} \rightarrow &  H^0(c^\ast (N_{Z/\mathbf P^n})) &&&&\\
&&&&\downarrow & &\downarrow&&&&\\
&&&& H^1(c^\ast (T_Y)) & &0 &&&&.
\end{array}\end{equation}

Let's define and analyze the diagram. 
Including all zero spaces, there are 5 rows and 6 columns. Second and third rows are parts of 
the long exact sequences from the short exact sequences 
of sheaves over $\mathbf P^1$, 

\begin{equation}\begin{array}{ccccccccc}
0 &\rightarrow & c^\ast(T_{X})&\rightarrow & c^\ast (T_{Y})&\rightarrow & c^\ast(N_{X/Y}) &\rightarrow & 0, \end{array}\end{equation}

\begin{equation}\begin{array}{ccccccccc}
0 &\rightarrow & c^\ast(T_{Z})&\rightarrow & c^\ast(T_{\mathbf P^n})&\rightarrow & c^\ast(N_{Z/\mathbf P^n}) &\rightarrow & 0. \end{array}\end{equation}

The third column is the part of the long sequence of the short exact sequence 
of sheaves over $\mathbf P^1$, 
\begin{equation}\begin{array}{ccccccccc}
0 &\rightarrow & c^\ast(T_{Y})&\rightarrow & c^\ast(T_{\mathbf P^n})&\rightarrow & c^\ast(N_{Y/\mathbf P^n}) &\rightarrow & 0. \end{array}\end{equation}
The isomorphism in the fourth column is from the adjuntion formula.  \bigskip

The  homomorphism $\mu_4$  is well-defined only if $c$  has unobstructed deformation on $Y$. Let's see this in the following. 
Because $c: \mathbf P^1\to Y$ is free, therefore  $c$  has unobstructed deformation on $Y$. So $H^1 (c^\ast(T_Y))=0$. 
Therefore the third column exact sequence in (4.9) splits, i.e. 
\begin{equation}
H^0(c^\ast(T_{\mathbf P^n}))\simeq  H^0(c^\ast (N_{Y/\mathbf P^n}))\oplus H^0( c^\ast (T_Y)).\end{equation}
We define $\mu_4$ to be the restriction of $\mu_3$ to the subspace, $H^0 (c^\ast(N_{Y/\mathbf P^n}))$  i.e.

\begin{equation}\begin{array}{ccc}

\mu_4: H^0 (c^\ast(N_{Y/\mathbf P^n})) &\rightarrow & H^0(c^\ast (N_{Z/\mathbf P^n}))\\
\alpha & \rightarrow & \alpha|_t+ c^\ast (T_{Z})|_t.
\end{array}\end{equation}
(But $\mu_4$ may not be zero.). 
In the sense of this splitting, the map $\mu_3$ also splits as  
\begin{equation}
\mu_3=\mu_1\oplus \mu_4.\end{equation}

Next we go further to use a construction to prove that $\mu_4$ is the zero map  due to the global generation.
This is just a specification of 
$$H^0 (c^\ast(N_{Y/\mathbf P^n}))$$ inside of 
$$H^0( c^\ast (T_Y)).$$
By (4.13), there are a Zariski open set $U\subset \mathbf P^1$ and a trivial subbundle
 $$E\subset  c^\ast(N_{Y/\mathbf P^n})$$ over $U$ such that
\begin{equation}
E\oplus  c^\ast(T_{Y})|_U= c^\ast(T_{\mathbf P^n})|_U
\end{equation} and
\begin{equation}
E\subset T_Z|_U.
\end{equation}

Let \begin{equation}
B=\{\sigma\in H^0(c^\ast(T_{\mathbf P^n})):\sigma|_U\in E\}.\end{equation}
So $B$ is a subspace of $H^0(c^\ast(T_{\mathbf P^n}))$. It generates a
sheaf $E_{\mathbf P^1}$ over $\mathbf P^1$ (in general $B$ could be zero, so did $E_{\mathbf P^1}$.). 
Because  $$c^\ast(T_{Y}),  c^\ast(T_{\mathbf P^n})$$ are generated by global section, 
the  sheaf $E_{\mathbf P^1}$ over $\mathbf P^1$ satisfies
\begin{equation} 
B\simeq  H^0(c^\ast (N_{Y/\mathbf P^n})).
\end{equation}
and (4.17) extends to 
\begin{equation}
E_{\mathbf P^1}\subset T_Z.
\end{equation}

We then define $\mu_4$ to be the restriction of $\mu_3$ to $B$
(this is the same as the definition before. But this time, the subspace  $H^0(c^\ast (N_{Y/\mathbf P^n}))$ is 
 uniquely identified.). By the condition (4.20), $\mu_4$ is the zero map. 

Thus  \begin{equation} 
Image(\mu_1)=Image (\mu_3).\end{equation}
Now we consider the third row, an exact sequence in (4.9). Since $Z$ is a generic hypersurface of $\mathbf P^n$, we can apply
 Theorem 1.1, [8] (or [7]) to obtain that  $H^1(c^\ast (T_{Z}))=0$. The exactness of the sequence implies that $\mu_3$ is surjective. 
So is  $\mu_1$.  Next we shift the focus to  the second row in (4.9).  We apply   $H^1(c^\ast (T_{Y}) ) = 0$ again to obtain that  
$H^1( c^\ast (T_{X}))=0$.  Using the exact sequence

\begin{equation}\begin{array}{ccccccccc}
0 &\rightarrow & T_{\mathbf P^1}&\rightarrow & c^\ast(T_{X})&\rightarrow & N_{c/X} &\rightarrow & 0. \end{array}\end{equation}
we obtain that \begin{equation}
H^1( c^\ast (T_{X}))=H^1(N_{c/X})=0.
\end{equation}

\end{proof}

\bigskip

\subsection{Free morphism}

In this subsection, we continue the existence result ---
the existence of irreducible rational curves of each degree $d$ on a generic
complete intersection Calabi-Yau 3-folds. 
From now on we change $Y$ to be a complete intersection as follows, 
\begin{equation}
Y=\cap_{i=1}^{r-1} div(g_i)
\end{equation}
for the fixed set of generic $g_i$.
 So $Y$,  which is equal to $X_4$,   is a smooth 4-dimensional complete intersection and Fano.
As before  for the fixed set of $g_i$, 
\begin{equation}
X_3=\cap_{i=1}^{r} div(g_i)
\end{equation}
 is a Calabi-Yau 3-fold contained in the Fano 4-fold $Y=X_4$. 
\bigskip

We would like to show a generic rational curve $c$ on 
$X_3$ is free on $Y$, in another word, such rational curves cover $Y$.   
 Thus we assume
\begin{ass}
 $c_0\in Hom_{bir}(\mathbf P^1, X_3)$ of $deg(c_0(\mathbf P^1))=d$ exists.
\end{ass}

  We'll work in a neighborhood of each component of the incidence schemes around $c_0$.
\bigskip

Let's resume all the notations in section 2. 
We let $g_1, \cdots, g_{r-1}$ be fixed (where $r=n-3$). Then the incidence scheme $I_{r-1}$ is also fixed. 
Let $\mathcal I_{r-1}$ be a component of $I_{r-1}$ containing $c_0$.
 Let $\Gamma$ be an irreducible component of the scheme
\begin{equation}
\{ (c, g_r)\in \mathcal I_{r-1}\times H^0(\mathcal O_{\mathbf P^n}(h_r)): c^\ast (g_r)=0\}
\end{equation}
containing $c_0$. 
Let \begin{equation}\begin{array}{ccc}
\pi: \mathcal I_{r-1} \times H^0(\mathcal O_{\mathbf P^n}(h_r)) &\rightarrow & M_d 
\end{array}\end{equation}
be the projection. 
Let $\mathcal I=\pi(\Gamma)$.  So $$\mathcal I\subset \mathcal I_{r-1}.$$

Let $0\in \mathbf P^1$ be a generic point. 
Let $\mathcal R$ be the closure of the open scheme, 
\begin{equation}
\mathring{\mathcal R}=\{ (c, y)\in \mathbf P(\mathcal I)\times Y: y=c(0)\}.
\end{equation}
The equation $y=c(0)$ means that $c$ is regular at $0$. 
\bigskip

We would like to show that
\bigskip

\begin{proposition}
The open scheme
 $\mathring{\mathcal R}$ is a rational map 
from $\mathbf P(\mathcal I)$ to $Y$ and dominates $Y$. 

\end{proposition}

\bigskip

\begin{proof}
Let $l$ be an integer satisfying
\begin{equation}
0\leq l\leq  r.
\end{equation}
Let $\mathcal I_l$ be a component of $ I_l$ containing $c_0$. 
Define 
\begin{equation}
\mathcal R_l\subset \mathbf P(\mathcal I_l)\times X_{n-l}\end{equation}
to be the closure of  
\begin{equation}
\mathring {R_l}=\{ (c, x): c(0)=x\}.
\end{equation}
for each $l$.

First we consider the case $l< r-1$.  Suppose $\mathcal R_l$ is onto $X_{n-l}$, where
$dim(X_{n-l})=n-l$.  Let $x\in X_{n-l}$ have coordinates $[0, \cdots, 0, 1]$.
Then  the dimension of the fibre $(\mathcal R_l)_x$ of $\mathcal R_l$ over the generic point $x\in X_{n-l}$ 
 satisfies 
\begin{equation}
dim( (\mathcal R_l)_x)\geq (h_{l+1}+\cdots+h_r)d.
\end{equation}

Hence after imposing $h_{l+1}d+1$ equations, 
$$g_{l+1}(c(t^{l+1}_1))=\cdots=g_{l+1}(c(t^{l+1}_{h_{l+1}d+1})=0,$$
we obtain that 
\begin{equation}
dim( (\mathcal R_{l+1})_x)\geq (h_{l+2}+\cdots+h_r)d-1 \geq 0.
\end{equation}
for all $x\in X_{n-l-1}$. 
Let $$Proj: \mathbf P(I_l)\to \mathbf P^n$$ be the projection.  Then the correspondence 
$\mathcal R_{l+1}$ sends $Proj(\mathcal R_{l+1})$ onto $X_{n-l-1}$ (not a map).  Next we show $ R_{l+1}$ is a rational map.

Applying Assumption 4.4, there is a $$c_0\in Proj(\mathcal R_{l+1})\cap M_{bir, d}.$$ 
Then generic $c\in  Proj (\mathcal R_{l+1})$
must be in $M_{bir, d}$ because $M_{bir, d}$ is an open set. 
Hence the correspondence $\mathcal R_{l+1}$ is a rational map.   
By the induction on the index $l$ for  $$l=0 \Rightarrow l=r-2, $$ we showed that
the irreducible component $$\mathring {\mathcal R}_{r-1}$$ dominates $Y$, through the rational evaluation map
$$c\to c(0).$$

\bigskip

To prove Proposition 4.5 it suffices to extend above proof to $l=r$, but the situation is slightly different. 

Let \begin{equation}
\mathbf t^r=(t^r_1, \cdots, t^r_{h_rd+1})\in sym^{h_rd+1}(\mathbf P^1)
\end{equation}
be generic. 
By the result of  above argument,  $\mathring {\mathcal R}_{r-1}$ dominates $Y$.
So for a generic $y\in Y$, 
\begin{equation}
dim ((\mathcal R_{r-1})_y)\geq h_rd.
\end{equation}

For the fixed $y$, we define $\mathcal S_y$ to be closure of 
\begin{equation}\begin{array}{c}
\mathring {\mathcal S_y}\subset \mathcal I_{r-1}\times \mathring{ H^0(\mathcal O_{\mathbf P^n}(h_r))}\\
\mathring{\mathcal S_y}=\{ (c, g_r): c\in I_{r-1}, g_r(c(t^r_1))=\cdots=g_r (c(t^r_{h_rd+1}))=0\}
\end{array}\end{equation}
where $\mathring{ H^0(\mathcal O_{\mathbf P^n}(h_r))}$ is an open set of space of hypersurface $g_r$ satisfying that
the Jacobian matrix of the functions in $c$ 
$$g_r(c(t^r_1)), \cdots, g_r (c(t^r_{h_rd+1}))$$
at $c\in I_{r-1}$ has full rank.  By the inequality (4.35), we count the dimension to obtain the 
$\mathcal S_y$ is non-empty. Hence 
\begin{equation}
\pi (\mathcal S_y)
\end{equation}
is non-empty for generic $y$.  This shows that the correspondence $\mathcal R$ is onto $Y$.
Now let's show it is a rational map. Notice a generic $c\in \pi (\mathcal S_y)$ is a generic element in 
\begin{equation}
\cup_{y\in Y}\pi (\mathcal S)=\mathcal I.
\end{equation}
Hence it suffices to show there is one regular point for $\mathcal R$. 
By Assumption 4.4, there is a $c_0\in {\mathcal S_y}$ such that $c_0$ is in $M_{bir, d}$, and $c_0(0)=y$ for some $y\in Y$ (not all $y$). 
Therefore $\mathcal R$ is regular at ONE point. With the same reason as above, $\mathcal R$ is a rational map.  
 Thus $\mathring {\mathcal R}$ dominates $Y$. 
This completes the proof.
\end{proof}
\bigskip

\bigskip

\begin{corollary}
Let $(c, y)\in \mathcal R$ be generic. Then $c$ is a free morphism on $Y$.

\end{corollary}
\bigskip

\begin{proof}
By Proposition 4.5, there is a Zariski open set $U$ of $\mathcal R$, such that the projection
$$Proj: U\to Proj(U)\subset Y$$ is smooth.
Hence its differential is onto. Let $(c, y)\in U$. Then the pull-back of the tangent bundle
\begin{equation}
c^\ast (T_Y) \simeq \oplus_k \mathcal O_{\mathbf P^1}(k)
\end{equation}
does not have negative summand, i.e. $k\geq 0$. Hence a  $c$ gives a free morphism
$\mathbf P^1\to Y$. 
\end{proof}
\bigskip

\subsection{Unobstructed deformation}
\bigskip

Notice that there is a degree $d$ rational curve $c\in Hom_{bir}(\mathbf P^1, X_3)$ for
 a generic $X_3$. 
By corollary 4.6,  it is a free morphism in $Y$, we apply Theorem 4.3.   Then the  second row of (4.9) implies a recursive formula
\begin{equation}
dim(H^0(c^\ast (X_{l+1})))=dim(H^0(c^\ast (X_{l})))-h_l d-1.
\end{equation}

Since $$dim(H^0(c^\ast (X_{0})))=(n+1)d+n, $$ 
\begin{equation} dim (H^0(c^\ast (X_3)))=(n+1)d+n-\sum_{k=1}^r h_k d-(n-3)=3.\end{equation}

Now we consider it from a different point of view.
Because $c$ is a birational map to its image,  there are finitely many points $t_i\in\mathbf P^1$ where the differential map

$$\begin{array}{ccc}
c_\ast: T_{t_i}\mathbf P^1 &\rightarrow & T_{c(t_i)}Y
\end{array}$$
is a zero map. Assume its vanishing order at $t_i$ is $m_i$ .  Let
\begin{equation}
m=\sum_i m_i.
\end{equation}
Let $s(t)\in H^0(\mathcal O_{\mathbf P^1}(m))$ such that $$div(s(t))=\Sigma_i m_it_i.$$

The sheaf morphism $c_\ast$ is injective and induces a composed morphism $\xi_s$ of sheaves
\begin{equation}\begin{array}{ccccc}
T_{\mathbf P^1} &\stackrel{c_\ast}  \rightarrow & c^\ast(T_{X_3}) & \stackrel{1\over s(t)}\rightarrow &
c^\ast(T_{X_3})\otimes \mathcal O_{\mathbf P^1}(-m).\end{array}
\end{equation}

It is easy to see that the induced bundle morphism $\xi_s$ is injective. Let 
\begin{equation}
N_m={c^\ast(T_{X_3})\otimes \mathcal O_{\mathbf P^1}(-m)\over \xi_s(T_{\mathbf P^1})}.
\end{equation}
Then
\begin{equation}
dim(H^0(N_m))=dim(H^0(c^\ast(T_{X_3})\otimes \mathcal O_{\mathbf P^1}(-m)))-3.
\end{equation}
On the other hand, three dimensional automorphism group of $\mathbf P^1$ gives a rise
to a 3-dimensional subspace  $K$ of  $$H^0(c^\ast(T_{X_3})).$$
By (4.41), $K=H^0(c^\ast(T_{X_3}))$. Over each point $t\in \mathbf P^1$, $K$ spans 
a one dimensional subspace. Hence 
\begin{equation} c^\ast(T_{X_3})
\simeq \mathcal O_{\mathbf P^1}(2)\oplus \mathcal O_{\mathbf P^1}(-k_1)\oplus \mathcal O_{\mathbf P^1}(-k_2),
\end{equation}
where $k_1, k_2$ are some positive integers.
This implies that
\begin{equation} dim \biggl(H^0(c^\ast T_{X_3}\otimes \mathcal O_{\mathbf P^1}(-m) )\biggr)=dim( H^0(\mathcal O_{\mathbf P^1}(2-m)).
\end{equation}
Then
\begin{equation} dim \biggl(H^0(c^\ast T_{X_3} \otimes \mathcal O_{\mathbf P^1}(-m))\biggr)=3-m.\end{equation}

Since $dim(H^0(N_m))\geq 0$, by the formula (4.45), $-m\geq 0$. 
By the definition of $m$, $m=0$.  Hence $c$ is an immersion. 
\par

Now $N_{c/X_3}$ is a vector bundle. By the the Calabi-Yau condition, 
\begin{equation}
N_{c/X_3}\simeq \mathcal O_{\mathbf P^1}(k)\oplus \mathcal O_{\mathbf P^1}(-2-k)
\end{equation}
where $k$ is non-positive. By Theorem 4.3, $H^1(N_{c/X_3})=0$. 
Hence $k=-1$.  Since $X_3$ is generic, the same proof is valid for all $c\in Hom_{bir}(\mathbf P^1, X_3)$.
Together with the existence in section 3, we complete the proof of Theorem 1.1.

\bigskip

\begin{appendices}

 \section{Vandermonde type of matrices}
 
 \bigskip
 
 In the appendix, we give proofs of lemmas dealing with Vandermonde type of matrices that are products of diagonal matrices and
Vandermonde matrices. The main purpose is to use rather elementary technique to
 glue
 matrices of smaller sizes, the Vandermonde type of matrices,  to obtain a matrix  of 
 larger size (which is a Jacobian matrix of a complete intersection).  
We assume that the linear algebra has  the ground  field $\mathbb C$. \bigskip

 First we introduce and study Vandermonde type of matrices, that are smaller matrices as entries in the large block matrices. 
  Let \begin{equation}
h \in H^0(\mathcal O_{\mathbf P^1}(v))
\end{equation}
be a polynomial, where $t\in \mathbf P^1$ is the variable, and
Let $$\mathbf t=(t_1, \cdots, t_{u})\in Sym^{u}(\mathbf P^1).$$
In the following in order to define matrices, we use the affine expression  in subsection 2.1. 
Let $\mathcal D_{u}$  be a diagonal matrix
\begin{equation}\mathcal D_u
=\left (\begin{array}{cccc} 
h(t_1) & 0 &\cdots & 0\\
0 & h(t_2)&  \cdots  &  0\\
\vdots &\cdots & \ddots & \vdots\\
 0 & \cdots & 0 & h(t_{u}) 
\end{array}\right)\end{equation}

\begin{definition}
\begin{equation}\mathcal V_0(h, \mathbf t, m)
=\mathcal D_{u}\left (\begin{array}{ccc} 
 t_1^{m} & \cdots  &  t_1\\
\vdots &\cdots & \vdots\\
t_{u}^{m} & \cdots & t_{u}
\end{array}\right).
\end{equation} 
 and 

\begin{equation}\mathcal V_{1}(h, \mathbf t, m)
=\mathcal D_{u}\left (\begin{array}{ccc} 
  t_1^{m} & \cdots  &  1\\
\vdots &\cdots & \vdots\\
t_{u}^m & \cdots & 1
\end{array}\right)\end{equation} 
We call them Vandermonde type of matrices of orders $0$ and $1$ respectively. 
\end{definition}

{\bf Remark}
 The Jacobian matrices $J_l$ of incidence schemes (we are interested in ) are all made of
Vandermonde type of matrices.
\bigskip

\begin{lemma}
If $h_1, h_2$ are two relatively prime polynomials with distinct zeros, and degrees $\geq d$,   then for generic 
$$(\mathbf t^1, \mathbf t^2)\in Sym^{d}(\mathbf P^1)\times Sym^{d+1}(\mathbf P^1)$$ 
the square matrix

\begin{equation} 
B=\left (\begin{array}{cc} 
 \mathcal V_0(h_1, \mathbf t^1, d)&\mathcal V_{1}(h_2, \mathbf t^1, d) \\
\mathcal V_{1}(h_1, \mathbf t^2, d)  & \mathcal V_{0}(h_2, \mathbf t^2, d) 
\end{array}\right)\end{equation}\\
is non-degenerate.

\end{lemma}
\bigskip

\begin{proof} It suffices to prove it for a special $\mathbf t^i$. 
Since $h_1$ has degree $\geq d$ and it is relatively prime to $h_2$, we can choose 
$ t^1_1, \cdots, t^1_d$ to be the zeros of $h_1(t)$ and $h_2(t)$ is non-zero at all points of
$\mathbf t^1, \mathbf t^2$. 
  Then 
\begin{equation}
B=\left (\begin{array}{cc} 
0&\mathcal V_{1}(h_2, \mathbf t^1, d) \\
\mathcal V_{1}(h_1, \mathbf t^2, d)  & \mathcal V_{0}(h_2, \mathbf t^2, d) 
\end{array}\right)\end{equation}

Then its determinant is
\begin{equation}
-|\mathcal V_{1}(h_2, \mathbf t^1) | |\mathcal V_{1}(h_1, \mathbf t^2) |.
\end{equation}
Since $\mathcal V_{1}(h_2, \mathbf t^1),  \mathcal V_{1}(h_1, \mathbf t^2)$ are types of Vandermonde matrices, the distinct
$t^i_j$ and nonvanishing $h_1, h_2$ imply the non-degeneracy. 
Therefore $B$ is non-degenerate. 

\end{proof}

\bigskip

\begin{lemma}
Let $h_3\in H^0(\mathcal O_{\mathbf P^1}(2d))$ and 
$h_1, h_2\in  H^0(\mathcal O_{\mathbf P^1}(d))$ be two vectors not on the same line through the origin. They are also pair wisely relatively prime,
and all zeros are distinct. Let $$ (\mathbf t^1, \mathbf t^2, \mathbf t^3)\in
 sym^{d}(\mathbf P^1)\times sym^{d}(\mathbf P^1)\times sym^{d+1}(\mathbf P^1)$$ be  generic.
Then the square matrix ( $(3d+1)\times (3d+1)$),
 \begin{equation}
B=\left (\begin{array}{ccc} 
\mathcal V_0(h_1, \mathbf t^1, d) &\mathcal V_0(h_2, \mathbf t^1, d)  & \mathcal V_1(h_3, \mathbf t^1, d)\\
\mathcal V_0(h_1, \mathbf t^2, d) &\mathcal V_0(h_2, \mathbf t^2, d)  & \mathcal V_1(h_3, \mathbf t^2, d)\\
\mathcal V_0(h_1, \mathbf t^3, d) &\mathcal V_0(h_2, \mathbf t^3, d)  & \mathcal V_1(h_3, \mathbf t^3, d)\\
\end{array}\right)\end{equation}
is non-degenerate.

\end{lemma}

\bigskip

\begin{proof} 
By the linear algebra in Proposition 2.10, [7],  we obtain that 
the matrix  $\mathcal V$, 
\begin{equation}
\mathcal V=\left (\begin{array}{cc} 
\mathcal V_0(h_1, \mathbf t^1, d) &\mathcal V_0(h_2, \mathbf t^1, d) \\
\mathcal V_0(h_1, \mathbf t^2, d) &\mathcal V_0(h_2, \mathbf t^2, d)  
\end{array}\right)\end{equation}
has full rank for any $(\mathbf t^1, \mathbf t^2)$ with $2d$ distinct points in $\mathbf P^1$. 

Then we rewrite $B$ as a block matrix, 
\begin{equation}
B=\left (\begin{array}{cc} 
\mathcal V  & \mathcal A\\
 \mathcal D & \mathcal V_1(h_3, \mathbf t^3, d) 
\end{array}\right)\end{equation}

where 
\begin{equation}
\mathcal D=\left (\begin{array}{cc} 
\mathcal V_0(h_1, \mathbf t^3, d) &\mathcal V_0(h_2, \mathbf t^3, d) 
\end{array}\right)\end{equation}
and \begin{equation}
\mathcal A=\left (\begin{array}{c} 
\mathcal V_1(h_3, \mathbf t^1, d) \\
\mathcal V_1(h_3, \mathbf t^2, d) 
\end{array}\right).\end{equation}

Then for generic $\mathbf t^1, \mathbf t^2, \mathbf t^3$,  $B$ is column-equivalent to

\begin{equation}
\left (\begin{array}{cc} 
\mathcal V  & 0 \\
 \mathcal D & -\mathcal D\mathcal V^{-1}\mathcal A+\mathcal V_1(h_3, \mathbf t^3, d)
\end{array}\right).\end{equation}
By choosing the points of $\mathbf t^1, \mathbf t^2$ to be the zeros of $h_3$ and $\mathcal V$ remains invertible, we obtain that
$$\mathcal A=0. $$
Thus we have a specialization 
$$ -\mathcal D\mathcal V^{-1}\mathcal A+\mathcal V_1(h_3, \mathbf t^3, d)=\mathcal V_1(h_3, \mathbf t^3, d)$$
which is non-degenerate. Thus for generic $\mathbf t^1, \mathbf t^2, \mathbf t^3$, 
$$ -\mathcal D\mathcal V^{-1}\mathcal A+\mathcal V_1(h_3, \mathbf t^3, d)$$
is non-degenerate, so is $B$. 

\end{proof}

\bigskip

We'll use two different types matrices  $B$ in Lemmas A. 2, A.3 as building blocks to obtain matrices of larger sizes.
We'll show that non-degeneracy of large matrices is the result of that  of the smaller, block matrices. 
\bigskip

\begin{lemma}
Let $B_{ij}, i=1, 2, j=1, \cdots, 6$ be $(3d+1)\times (d+1)$ matrices. 
Assume
\par
(1) \begin{equation}\left (\begin{array}{cc} 
B_{13} &B_{14} 
\end{array}\right)
=\left (\begin{array}{cc} 
B_{23} &B_{24} 
\end{array}\right)\end{equation}

(2) \begin{equation}\left (\begin{array}{ccc} 
B_{13} &B_{14} & B_{16}
\end{array}\right)\end{equation} and 
\begin{equation}\left (\begin{array}{ccc} 
B_{11} &  B_{12}-B_{22} &  B_{25} 
\end{array}\right)\end{equation}
have full rank, \par

(3) the columns of \begin{equation}\left (\begin{array}{cc} 
B_{13} &B_{14} 
\end{array}\right)\end{equation} 
 span the columns of 
\begin{equation}\left (\begin{array}{cc} 
B_{11} &B_{12} 
\end{array}\right)\end{equation} 

\par
Then 

\begin{equation}J=\left (\begin{array}{cccccc} 
B_{11} &  B_{12} &  B_{13} & B_{14} &0 & B_{16}\\
0 & B_{22} & B_{23} & B_{24} & B_{25} & B_{26}
\end{array}\right)\end{equation}

is non-degenerate.

\end{lemma}

\bigskip

\begin{proof} In the following, we apply the column and row operations to the matrix. They will not change its rank.
Applying the row operations on the matrix $J$ we obtain it is row-equivalent to
\begin{equation}J=\left (\begin{array}{cccccc} 
B_{11} &  B_{12} &  B_{13} & B_{14} &0 & B_{16}\\
-B_{11} & B_{22}-B_{12} & 0 & 0 & B_{25} & B_{26}-B_{16}
\end{array}\right)\end{equation}
By the condition (3), it is column-equivalent to
\begin{equation}J=\left (\begin{array}{cccccc} 
0 &  0 &  B_{13} & B_{14} &0 & B_{16}\\
-B_{11} & B_{22}-B_{12} & 0 & 0 & B_{25} & B_{26}-B_{16}
\end{array}\right)\end{equation}
By the condition (2), it is further column equivalent to 
\begin{equation}J=\left (\begin{array}{cccccc} 
0 &  0 &  B_{13} & B_{14} &0 & B_{16}\\
-B_{11} & B_{22}-B_{12} & 0 & 0 & B_{25} & 0 
\end{array}\right)\end{equation}
By the condition (2), we obtain it has full rank.

\end{proof}

\bigskip

\begin{lemma}
Consider the block matrix
\begin{equation}J=\left (\begin{array}{ccccccc} 
B_{11} &  B_{12} &  B_{13} & B_{14} & B_{15} & 0 & 0\\
0 &  0 & B_{23} & B_{24} & B_{25} & B_{26} & 0\\
0 &  0 & 0 & B_{34} & B_{35} & B_{36} & B_{37}
\end{array}\right)\end{equation}
where the entries in first row are $(3d+1)\times (d+1)$ matrices and all the rest
are $(2d+1)\times (d+1)$ matrices. We assume
\par
(1) \begin{equation}\left (\begin{array}{ccc} 
B_{24} &B_{25}  & B_{26}
\end{array}\right)=\left (\begin{array}{ccc} 
B_{34} &B_{35} &B_{36} 
\end{array}\right).
\end{equation} 
 (2) \begin{equation}\left (\begin{array}{ccc} 
B_{11} &B_{12} & B_{15}
\end{array}\right)\end{equation},  
\begin{equation}\left (\begin{array}{cc} 
B_{23} &  -B_{37}  
\end{array}\right)\end{equation} and

\begin{equation}\left (\begin{array}{cc} 
B_{34} & B_{36}  
\end{array}\right)\end{equation} 

have full rank, \par

(3) the columns of \begin{equation}\left (\begin{array}{cc} 
B_{11} &B_{12} 
\end{array}\right)\end{equation} 
 span the columns of 
\begin{equation}\left (\begin{array}{cc} 
B_{13} &B_{14} 
\end{array}\right),\end{equation} 
then the matrix $J$ has full rank.
\end{lemma}

\bigskip

\begin{proof} The proof uses row and column operations which are  the same as those in Lemma A.4.
So in the following we just list the equivalent matrices in the reduction.
\begin{equation}J=\left (\begin{array}{ccccccc} 
B_{11} &  B_{12} &  B_{13} & B_{14} & B_{15} & 0 & 0\\
0 &  0 & B_{23} & B_{24} & B_{25} & B_{26} & 0\\
0 & 0 & 0 & B_{34} & B_{35} & B_{36} & B_{37}
\end{array}\right)\end{equation}
$$\Downarrow\scriptstyle{row \ operations} $$
\begin{equation}\left (\begin{array}{ccccccc} 
B_{11} &  B_{12} &  B_{13} & B_{14} & B_{15} & 0 & 0\\
0 &  0 & B_{23} & 0 & 0& 0 & -B_{37}\\
0 & 0 & 0 & B_{34} & B_{35} & B_{36} & B_{37}
\end{array}\right)\end{equation}
$$\Downarrow{column\ operations}$$

\begin{equation}J=\left (\begin{array}{ccccccc} 
B_{11} &  B_{12} &  0 & 0 & B_{15} & 0 & 0\\
0 &  0 & B_{23} & 0 & 0 & 0 & -B_{37}\\
0 & 0 & 0 & B_{34} & B_{35} & B_{36} & B_{37}
\end{array}\right)\end{equation}
$$\Downarrow{column\ operations}$$
\begin{equation}J=\left (\begin{array}{ccccccc} 
B_{11} &  B_{12} &  0 & 0 & B_{15} & 0 & 0\\
0 &  0 & B_{23} & 0 & 0 & 0 & -B_{37}\\
0 & 0 & 0 & B_{34} & 0 & B_{36} & 0
\end{array}\right).\end{equation}
The last matrix has full rank. 

\end{proof}

\begin{lemma}

Consider the block matrix
\begin{equation}J=\left (\begin{array}{cccc} 
B_{11} &B_{12} &0 & 0 \\
0 &B_{22} & B_{23} &0  \\
0 &0& B_{33} & B_{34}\\
0 &B_{42} &0  & B_{44}\\
\end{array}\right),\end{equation} 
where $B_{ij}$ are non degenerate square matrices of the same size. If
\begin{equation}\left (\begin{array}{cc} 
-B_{22} B_{23}^{-1} B_{33}  &B_{34} \\
B_{42} & B_{44}
\end{array}\right)\end{equation} 
is non-degenerate, so is $J$.
\end{lemma}
\bigskip

\begin{proof} 
 As before we perform row and column operations  on $J$ to obtain 
\begin{equation}\left (\begin{array}{cccc} 
B_{11} &0 &0 & 0 \\
0 &B_{22} & B_{23} &0  \\
0 &0 & B_{33} & B_{34}\\
0 &B_{42} &0  & B_{44}
\end{array}\right)\end{equation} 
$$\Downarrow$$

\begin{equation}\left (\begin{array}{cccc} 
B_{11} &0 &0 & 0 \\
0 &0 & B_{23} &0  \\
0 & -B_{22} B_{23}^{-1} B_{33} & B_{33} & B_{34}\\
0 &B_{42} &0  & B_{44}
\end{array}\right)\end{equation} 
$$\Downarrow$$ 
\begin{equation}\left (\begin{array}{cccc} 
B_{11} &0 &0 & 0 \\
0 &0 & B_{23} &0  \\
0 & -B_{22} B_{23}^{-1} B_{33} & 0 & B_{34}\\
0 &B_{42} &0  & B_{44}
\end{array}\right)\end{equation} 

By the assumption, \begin{equation}\left (\begin{array}{cc} 
-B_{22} B_{23}^{-1} B_{33}  &B_{34} \\
B_{42} & B_{44}
\end{array}\right)\end{equation} 
is non-degenerate, so is $J$.

\end{proof} 
\bigskip

\end{appendices}

\end{document}